\documentclass[11pt,twoside]{article}
\usepackage{mathbbold}
\usepackage{amsmath, amssymb, mathrsfs, graphicx}

\def\scr{\mathscr}

\def\az{\alpha}  \def\bz{\beta}
    \def\dz{\delta}
    \def\fz{\varphi}
\def\gz{\gamma}  \def\kz{\kappa}
\def\lz{\lambda}

        \def\uz{\theta}
\def\vz{\varepsilon}

\def\ffz{\Phi}    
  \def\ooz{\Omega}

\def\qd{\quad}
\def\qqd{\qquad}

\newcommand{\mathsym}[1]{{}}

\def\scr{\mathscr}

\def\le{\leqslant}
\def\ge{\geqslant}

\font\cms=cmss9 scaled \magstep1

\def\nnd{\noindent}

\def\thm{\nnd\bg{thm1}}
\def\crl{\nnd\bg{crl1}}

\def\prp{\nnd\bg{prp1}}

\def\xmp{\nnd\bg{xmp1}}

\def\rmk{\nnd\bg{rmk1}}

\def\dethm{\end{thm1}}
\def\decrl{\end{crl1}}

\def\deprp{\end{prp1}}
\def\dexmp{\end{xmp1}}

\def\dermk{\end{rmk1}}

\def\prf{\medskip \noindent {\bf Proof}. }
\def\deprf{\quad $\square$ \medskip}
\def\bg{\begin}
\def\be{\bg{equation}}
\def\de{\end{equation}}

\def\dear{\end{eqnarray}}
\def\lb{\label}
\def\ct{\cite}

\newcommand{\rf}[2]{[\ref{#1}; #2]}

\def\den{\end{enumerate}}

\def\d{\text{\rm d}}

\begin{document}

\allowdisplaybreaks[4]
\thispagestyle{empty}
\renewcommand{\thefootnote}{\fnsymbol{footnote}}

\noindent {Front. Math. China}

\vspace*{.5in}
\begin{center}
{\bf\Large Lower bounds of principal eigenvalue in dimension one}
\vskip.15in {Mu-Fa Chen}
\end{center}
\begin{center} (Beijing Normal University, Beijing 100875, China) \end{center}
\vskip.1in

\markboth{\sc Mu-Fa Chen}{\sc Lower bound of the principal eigenvalue}


\date{}


\footnotetext{Received November 30, 2011; accepted May 31, 2012}

\footnotetext{{\it Key words and phases}.\quad
Principal eigenvalue, lower estimate, variational formula, one-dimensional diffusion,
birth--death process.}
\footnotetext{Research supported in part by the
         National Natural Science Foundation of China (No. 11131003),
         and by the ``985'' project from the Ministry of Education in China.
}

\bigskip

\begin{abstract}
For the principal
eigenvalue with bilateral Dirichlet boundary condition, the so-called basic
estimates were originally obtained by capacitary method.
The Neumann case (i.e., the ergodic case) is even harder,
and was deduced from the Dirichlet one plus a use of duality and the coupling method.
In this paper, an alternative and more direct proof for the basic estimates
is presented. The estimates in the Dirichlet case are then improved by a typical application of
a recent variational formula. As a dual of the Dirichlet case, the refine problem for bilateral Neumann boundary condition is also treated. The paper starts with
the continuous case (one-dimensional diffusions) and ends at the discrete one
(birth--death processes). Possible generalization of the results
studied here is discussed at the end of the paper.
\end{abstract}

\medskip

\section{Introduction (continuous case)}
Consider an elliptic operator
\be L=a(x)\frac{\d^2}{\d x^2}+ b(x)\frac{\d}{\d x} \de
(with $a>0$) on $E:=(-M, N)$\, $(M, N\le \infty)$. Define a function $C(x)$:
$$C(x)=\!\int_{o}^x \frac b a,\qqd x\in E,$$
where $o\in E$ is a reference point. Here and in what follows, the
Lebesgue measure $\d x$ is often omitted. It is convenient for us to
define two measures $\mu$ and $\nu$:
\be \mu(\d x)=\frac{e^{C(x)}}{a(x)}\d x,\qqd
\nu(\d x)= e^{-C(x)}\d x. \de
As usual, the norm on $L^2(\mu)$ is denoted by $\|\cdot\|$. Define
$$\aligned
 {\scr A}(-M, N)&= \text{the
set of absolutely continuous functions
on $(-M, N)$},\\
{\scr A}_0(-M, N)&= \{f\in {\scr A}(-M, N): f\text{ has a compact
support} \},\\
D(f)&=\int_{-M}^N {f'}^2 e^C,\qqd f\in {\scr A}(-M, N),\;\;  M, N\le \infty.
\endaligned
$$
Here $D(f)$ is allowed to be $\infty$.
We are interested in the following eigenvalues:
\begin{align}
\lz^{\text{\rm DD}}&=\inf\{D(f): f\in {\scr A}_0(-M, N), \; \|f\|=1\},\\
\lz^{\text{\rm NN}}&=\inf\{D(f): f\in {\scr A}(-M, N), \;\mu(f)=0,\;
\|f\|=1\},
\end{align}
where  $\mu(f)=\int_E f \d \mu$.

The basic estimates, of $\lz^{\text{\rm DD}}$ for instance, given in \ct{cmf10} are as follows:
\be\big(4\,\kz^{\text{\rm DD}}\big)^{-1}\le \lz^{\text{\rm DD}}\le \big(\kz^{\text{\rm DD}}\big)^{-1},\lb{40}\de
where
\be\big(\kz^{\text{\rm DD}}\big)^{-1}\!\!= \inf_{x<y}
\big[\nu (-M, x)^{-1}\! + \nu (y, N)^{-1}\big]\, \mu (x, y)^{-1},\;\;
\mbox{$\mu(x, y)\!:=\!\!\int_x^y \d\mu$}.\lb{05}\de
The proof for the upper estimate is already straightforward, simply using
the classical variational formula for $\lz^{\text{\rm DD}}$ (cf. \rf{cmf10}{Proof (b) of Theorem 8.2}).
However, the proof for the lower estimate is much harder and deeper, using
capacity theory (cf. \rf{cmf10}{Sections 8, 10}). Even through the capacitary tool is suitable
in a general setup (cf. \ct{fmut}, \rf{cmf05}{Theorems 7.1 and 7.2},
\rf{mav}{Chapter 2}), it is still expected to have a direct proof (avoiding capacity)
in such a concrete situation. This is done at the beginning of the next section.
Surprisingly, the simple proof also works in the ergodic case for which the original proof
is based on (\ref{40}) plus a use of the duality and the coupling technique.
The main body of the paper is devoted to an improvement of the basic lower estimate given in (\ref{40}),
as stated in Corollary \ref{t11} below. The result can be regarded as a typical application  of a recent variational formula
(\rf{cmf11}{Theorem 4.2} or Theorem \ref{t21} below).
This note is an addition to the recent papers \ct{cmf10, cmf11} from which one can find
the motivation of the study on the topic and further references. It is remarkable that the new result makes
the whole analytic proof for the basic estimates more elementary.

Here is our first main result which is a refinement of \rf{cmf11}{Corollary 4.3}.

\crl\lb{t11} {\cms
\begin{itemize}
\item [(1)] We have
$$\lz^{\text{\rm DD}}\ge \big({\underline\kz}^{\text{\rm DD}}\big)^{-1}\ge \big(4\,\kz^{\text{\rm DD}}\big)^{-1},$$
where $\kz^{\text{\rm DD}}$ is given in (\ref{05}) and ${\underline\kz}^{\text{\rm DD}}$ is defined by (\ref{11}) below.
\item [(2)] Let $\mu(-M, N)<\infty$. Then assertion (1) holds if the codes DD are replaced by NN
\big(for instance, $\lz^{\text{\rm NN}}\ge \big({\underline\kz}^{\text{\rm NN}}\big)^{-1}$\big)
and the measures $\mu$ and $\nu$ are exchanged.
\end{itemize}}
\decrl

The remainder of the paper is organized as follows. In the next section,
we present shortly an alternative proof of the estimates in (\ref{40}).
The proof shows one of the main new ideas of the paper. Then we prove Corollary \ref{t11}.
Two illustrating examples are also included in this section.
The discrete analog
of Corollary \ref{t11} is presented in the third section.

\section{Proofs and Examples}

\medskip \noindent {\bf Proof of (\ref{40})}.

Let $\uz\in (-M, N)$ be a reference point. Define
$$\dz_{\uz}^-\!=\!\!\sup_{z\in (-M,\, \uz)}\nu(-M,\, z)\, \mu(z,\, \uz),\qqd
\dz_{\uz}^+\!=\!\!\sup_{z\in (\uz,\, N)}  \mu(\uz,\, z)\,\nu(z,\, N).$$
As will be remarked in the next section, we may assume that
$\dz_{\uz}^{\pm}<\infty$. Otherwise, the problem becomes
either trivial or degenerated. Next, denote by $\lz_{\uz}^{\pm}$ the principal
eigenvalue on $(-M, \uz)$ and $(\uz, N)$, respectively, with common reflecting (Neumann)
boundary at $\uz$ and absorbing (Dirichlet) boundary at $-M$ (and $N$) provided $M< \infty$
($N<\infty$). Actually, by an approximating procedure, one may assume that $M, N<\infty$
(cf. \rf{cmf10}{Proof of Corollary 7.9}).
Next, by a splitting technique, one may choose $\uz=\bar\uz$ to be the unique solution to the
equation $\lz_{\uz}^-=\lz_{\uz}^+$. Then they coincide with $\lz^{\text{\rm DD}}$ since by
\rf{czz03}{Theorem 1.1}, we have
$$\lz_{\uz}^-\wedge \lz_{\uz}^+\le \lz^{\text{\rm DD}}\le \lz_{\uz}^-\vee\lz_{\uz}^+$$
for every $\uz\in (-M, N)$,
where $x\wedge y =\min\{x, y\}$ and dually $x\vee y =\max\{x, y\}$.
Alternatively, $\bar\uz$ is the root of the derivative of the eigenfunction of $\lz^{\text{\rm DD}}$
by \rf{czz03}{Proposition 1.3} and the monotonicity of the eigenfunctions of $\lz_{\bar\uz}^{\pm}$. From now on in this proof, we fix this $\bar\uz$.
For given $\vz>0$, let $\bar x<\bar\uz$ and $\bar y>\bar\uz$ satisfy
$$\nu(-M,\,\bar x)\,\mu(\bar x,\, \bar\uz)\ge \dz_{\bar\uz}^--\vz,\qqd
 \mu(\bar\uz,\, \bar y)\, \nu(\bar y,\, N) \ge \dz_{\bar\uz}^+-\vz,$$
respectively. As a continuous analog of \rf{cmf00}{Theorem 1.1}, we have
$$\big[\big(\lz^{\text{\rm DD}}\big)^{-1}
=\big] \qqd \big(\lz_{\bar\uz}^{+}\big)^{-1}\le 4\,\dz_{\bar\uz}^+\qqd
\big[\le 4  \mu(\bar\uz,\, \bar y)\, \nu(\bar y,\, N)+ 4\vz\big].$$
Hence,
$$\big[\big(\lz^{\text{\rm DD}}\big)^{-1}-4\vz\big]\nu(\bar y,\, N)^{-1}\le 4 \mu(\bar\uz,\, \bar y).$$
In parallel, we have
$$\big[\big(\lz^{\text{\rm DD}}\big)^{-1}-4\vz\big]\nu(-M,\, \bar x)^{-1}\le 4 \mu(\bar x,\,\bar\uz).$$
Summing up the last two inequalities, it follows that
$$\big[\big(\lz^{\text{\rm DD}}\big)^{-1}-4\vz\big]\big[ \nu(-M,\, \bar x)^{-1}
+\nu(\bar y,\, N)^{-1}\big]
\le 4 \mu(\bar x,\,\bar y).$$
That is,
$$\big(\lz^{\text{\rm DD}}\big)^{-1}-4\vz
\le 4 \big[ \nu(-M,\, \bar x)^{-1}
+\nu(\bar y,\, N)^{-1}\big]^{-1}
\mu(\bar x,\,\bar y).$$
In view of (\ref{05}), the right-hand side is bounded from above by $4 {\kz}^{\text{\rm DD}}.$
Since $\vz$ is arbitrary, we have proved the lower estimate in (\ref{40}).
A direct proof for the upper one in (\ref{40}) is presented in \rf{cmf10}{Proof (b) of Theorem 8.2}.
\deprf

\medskip \noindent {\bf Proof of the dual of (\ref{40})}:
$$\big(4\,\kz^{\text{\rm NN}}\big)^{-1}\le \lz^{\text{\rm NN}}\le \big(\kz^{\text{\rm NN}}\big)^{-1},$$
where
$$\big(\kz^{\text{\rm NN}}\big)^{-1}\!\!= \inf_{x<y}
\big[\mu (-M, x)^{-1}\! + \mu (y, N)^{-1}\big]\, \nu (x, y)^{-1}.$$
By exchanging ``Neumann'' and ``Dirichlet'', the splitting
point $\uz=\bar\uz$ is now a common Dirichlet boundary
and $-M$ becomes Neumann boundary if $M<\infty$ (and so is $N$). In other words,
$\bar\uz$ is the unique root of the eigenfunction of $\lz^{\text{\rm NN}}$.
Now, in the proof above, we need only to use \rf{cmf00}{Theorem 3.3}
instead of \rf{czz03}{Theorem 1.1} and making the exchange of $\mu$ and $\nu$.
We have thus returned to the role mentioned in \ct{cmf11}: exchanging the boundary condition
``Neumann'' and ``Dirichlet'' simultaneously
leads to the exchange of the measures $\mu$ and $\nu$.

Here is a direct proof for the upper estimate. Given $x, y\in (-M, N)$ with $x<y$, let $\bar\uz=\bar\uz (x, y)$ be the unique solution to the equation
$$\aligned
&\mu(-M, x)\nu(x, \uz)+\int_x^{\uz}\mu(\d z)\nu(z, \uz)\\
&\qd =\mu(y, N)\nu(\uz, y)+\int_{\uz}^y\mu(\d z)\nu(\uz, z),\qqd \uz\in (x, y).\endaligned$$
Next, define
$$f(z)=-\mathbbold{1}_{\{z\le \bar\uz\}}\nu\big(x\vee z, \bar\uz\big)
+\mathbbold{1}_{\{z>\bar\uz\}}\nu\big(\bar\uz, y\wedge z\big).$$
Then $\mu(f)=0$ by the definition of $\bar\uz$. We have
$$\int_{-M}^N \big|f'\big|^2 e^C
=\nu\big(x, \bar\uz\big)+\nu\big(\bar\uz, y\big)
=\nu\big(x, y\big).$$
Moreover,
$$\aligned
\int_{-M}^N \big(f-\pi(f)\big)^2\d\mu&= \int_{-M}^N f^2\d\mu\\
&> \int_{-M}^x f^2\d\mu+\int_{y}^N f^2\d\mu\\
&=\mu(-M, x)\, \nu\big(x, \bar\uz\big)^2+\mu(y, N)\,\nu\big(\bar\uz, y\big)^2.
\endaligned$$
Note that the function
$$\gz(x)=\az x^2 +\bz (1-x)^2, \qqd x\in (0, 1),\;
\az, \bz>0$$
achieves its minimum $\big(\az^{-1}+\bz^{-1}\big)^{-1}$
at $x^*=(1+\bz/\az)^{-1}$. As an application of this result with
$$\az=\mu(-M, x),\qd \bz=\mu(y, N),\qd x=\nu\big(x, \bar\uz\big)/\nu\big(x, y\big),$$
we get
$$\int_{-M}^N \big(f-\pi(f)\big)^2\d\mu
\ge \frac{\nu(x, y)^2}{\mu(-M, x)^{-1}+\mu(y, N)^{-1}}.$$
Hence
$$\frac{\int_{-M}^N \big(f-\pi(f)\big)^2\d\mu}{\int_{-M}^N \big|f'\big|^2 e^C}
\ge \frac{\nu(x, y)}{\mu(-M, x)^{-1}+\mu(y, N)^{-1}}.$$
Making supremum with respect to $x<y$, we obtain the required
$\kz^{\text{\rm NN}}$.
\deprf

It is remarkable that although the last proof is in parallel to the previous one,
it does not depend on (\ref{40}). This is rather lucky since in other cases,
part (2) of Corollary \ref{t11} for instance, we do not have such a direct
proof.

From now on, unless otherwise stated, we restrict ourselves to the Dirichlet case.
For fixed $\uz$, much knowledge on $\lz_{\uz}^{\pm}$ is known (variational formulas,
approximating procedure and so on, refer to \ct{cmf05, cmf10} for instance).
Of which, only a little is used in the proof above.
For instance, by \rf{czz03}{Corollary 1.5}, we have
$$\Big(\sup_{\uz}\big[\dz_{\uz}^-\wedge \dz_{\uz}^+\big]\Big)^{-1}
\ge \lz^{\text{\rm DD}}
\ge \Big(4\inf_{\uz}\big[\dz_{\uz}^-\vee \dz_{\uz}^+\big]\Big)^{-1}.$$
Thus, if we choose $\bar\uz$ to be the solution of
equation $\dz_{\uz}^-= \dz_{\uz}^+$, then we obtain
$$\big(\dz_{\bar\uz}^-\big)^{-1}
\ge \lz^{\text{\rm DD}}
\ge \big(4\dz_{\bar\uz}^-\big)^{-1}$$
which is even more compact than (\ref{40}) in view of the comparison of $\kz^{\text{\rm DD}}$ and $\dz_{\uz}^{\pm}$.
The problem is that $\bar\uz$, especially the one used in the first proof of this section, is usually not
explicitly known and so a large part of the known results for $\lz_{\bar\uz}^{\pm}$ are not practical.
To overcome this difficulty, the first proof above uses two parameters $x$ and $y$
to get $\kz^{\text{\rm DD}}$ and then to obtain the explicit lower
estimates (\ref{40}). For our main result Corollary {\ref{t11}}, the fixed
point $\bar\uz$ used in the proof of (\ref{40}) is replaced by its mimic given in (\ref{08}) below for suitable
test function $f$. The difference is that equation (\ref{08}) is explicit but
not the one for $\bar\uz$ used in the first proof above.

\medskip \noindent {\bf Proof of Corollary \ref{t11}\,(1)}.

By \ct{cmf10} or \ct{cmf11}, we have known that part (2) of Corollary \ref{t11}
is a dual of part (1). Hence in what follows, we need study part (1) only.

The first inequality in part (1) comes from \rf{cmf11}{Corollary 4.3}.
Thus, it suffices to prove the last inequality in part (1).

Even though it is not completely necessary,
we assume that $M, N<\infty$ until the last paragraph of the proof.

For a given $f\in {\scr C}_+:$
$${\scr C}_+\!=\!\{f\!\in\! {\scr C}(-M, N)\!: f\!>\!0 \;\text{on}\; (-M, N),\; f(-M+0)\!=\!0 \text{ and } f(N-0)\!=\!0\},$$
define
\begin{align}
h^-(z)&\!=\!h_f^-(z)
\!=\!\!\int_{-M}^{z}\! e^{-C(u)}\d u \int_u^{\uz}\!\frac{e^C f}{a},\qqd
 z\le \uz, \\
h^+(z)&\!=\!h_f^+(z)
\!=\!\!\int_{z}^{N}\! e^{-C(u)}\d u \int_{\uz}^u\!\frac{e^C f}{a},\qqd z>\uz, \lb{8}
\end{align}
where $\uz=\uz (f)\in (-M, N)$ is the unique root of the equation:
\be h^-(\uz)= h^+(\uz)\lb{08}\de
provided $h_f^{\pm}<\infty$. The uniqueness of $\uz$ should be clear since on $(-M, N)$,
as a function of $\uz$, $h^-(\uz)$ is continuously increasing from zero to $h^-(N-0)>0$
and $h^+(\uz)$ is continuously decreasing from $h^+(-M+0)>0$ to zero.
Next, define
$$I\!I^{\pm}(f)= h^{\pm}/f.$$
Then we have the following variational formula.
\thm\lb{t21}{\bf\rf{cmf11}{Theorem 4.2\,(1)}}\;\;{\cms  Assume that $\nu(-M, N)<\infty$. Then
\be\lz^{\text{\rm DD}}= \sup_{f\in {\scr C}_+}\Big\{\Big[\inf_{z\in (-m,\,\uz)}I\!I^-(f)(z)^{-1}\Big]
\bigwedge \Big[\inf_{z\in (\uz, N)}I\!I^+(f)(z)^{-1}\Big]\Big\}.\lb{9}\de
}\dethm

We remark that in the original statement of \rf{cmf11}{Theorem 4.2\,(1)}, the boundary
condition ``$f(-M+0)\!=\!0 \text{ and } f(N-0)\!=\!0$'' is ignored. The condition is added
here for the use of the operators $I^{\pm}$ (different from $I\!I^{\pm}$) to be defined later.
However, the conclusion (\ref{9}) remains true since the eigenfunction of $\lz^{\text{\rm DD}}$
does satisfy this condition.

We now fix $x<y$ and let $f=f^{x, y}$:
\be f^{x, y}(s)=
\begin{cases}
\sqrt{\fz^+(y)\fz^-(s\wedge x)/\fz^-(x)},\qd & s\le y\\
\sqrt{\fz^+(s)}, & s\ge y,
\end{cases}
\de
where
$$\fz^-(s)=\nu (-M, s)\qd\text{ and }\qd\fz^+(s)=\nu(s, N).$$
Certainly, here we assume that $\fz^{\pm}<\infty$ (which is automatic whenever $M, N< \infty$).
Clearly, $f^{x, y}\in {\scr C}_+$. Here we are mainly interested in those pair $\{x, y\}$ having the property
$x< \uz< y$.  As proved in \ct{cmf11}, the quantity ${\underline\kz}^{\text{\rm DD}}$:
\begin{align}
{\underline\kz}^{\text{\rm DD}}
&=\inf_{x<y}\Big[\sup_{z\in (-M,\, \uz)}I\!I^-(f^{x, y})(z)\Big]
\bigvee \Big[\sup_{z\in (\uz, N)}I\!I^+(f^{x, y})(z)\Big]\lb{11}
\end{align}
used in Corollary \ref{t11}\,(1) has an explicit expression:
\begin{align}
&\inf_{x< y}
\bigg\{
\sup_{z\in (-M,\, x)}\bigg[\frac{1}{\sqrt{\fz^-(z)}}\,
\mu\Big((\fz^-)^{3/2} \mathbbold{1}_{(-M,\, z)}\Big)
+ \sqrt{{\fz^-(z)}}\,
\mu\Big(\sqrt{\fz^-}\, \mathbbold{1}_{(z,\, x)}\Big)\nonumber\\
&\qqd\qqd\qqd\qd +\sqrt{\fz^-(z)\fz^-(x)}\,\mu(x,\, \uz)\bigg]\nonumber\\
&\qqd\qqd\bigvee \bigg[\frac{1}{\sqrt{\fz^-(x)}}\,
\mu\Big((\fz^-)^{3/2} \mathbbold{1}_{(-M,\, x)}\Big)
+ \mu\Big(\fz^-\, \mathbbold{1}_{(x,\, \uz)}\Big)\bigg]\nonumber\\
&\qqd\qqd \bigvee\sup_{z\in (y,\, N)} \bigg[\frac{1}{\sqrt{\fz^+(z)}}\mu\Big((\fz^+)^{3/2} \mathbbold{1}_{(z,\, N)}\Big)
+ \sqrt{\fz^+(z)}\,\mu\Big(\sqrt{\fz^+}\, \mathbbold{1}_{(y,\, z)}\Big)\nonumber\\
&\qqd\qqd\qqd\qd\qqd +\sqrt{\fz^+(z)\fz^+(y)}\,\mu(\uz,\, y)\bigg]
\bigg\}.\nonumber
\end{align}
We have thus sketched the original attempt (cf. \rf{cmf11}{Corollary 4.3}) to prove Corollary \ref{t11}\,(1). The study was stopped here since we were unable
to compare this long expression with $4\,\kz^{\text{\rm DD}}$.

Before moving further, let us make a remark on (\ref{08}). As proved in \rf{cmf11}{(31)}, for fixed $x$ and $y$, equation (\ref{08}) is
equivalent to the following one.
\begin{align}
\!\!\!\!\frac{\mu\big((\fz^-)^{3/2} \mathbbold{1}_{(-M,\, x)}\big)}{\sqrt{\fz^-(x)}}\,
\!+\! \mu\big(\fz^-\, \mathbbold{1}_{(x,\, \uz)}\big)
\!=\!\frac{\mu\big((\fz^+)^{3/2} \mathbbold{1}_{(y,\, N)}\big)}{\sqrt{\fz^+(y)}}\,
\!+\! \mu\big(\fz^+\, \mathbbold{1}_{(\uz,\, y)}\big)\lb{12}.\end{align}
The quantity in (\ref{12}) is actually the ratio
$$\frac{h^-(\uz)}{f^{x, y}(x)}=\frac{h^+(\uz)}{f^{x, y}(y)}$$
(cf. \rf{cmf11}{(34)}) noting that $f^{x, y}$ is a constant on $[x, y]$:
$$f^{x, y}(x)=f^{x, y}(y)=\sqrt{\fz^+(y)}.$$
Next, note that the left-hand and the right-hand sides of (\ref{12}) are monotone, with respect to $x$ and $y$ respectively, since each of their derivatives does not change its sign:
$$-\frac{e^{-C(x)}}{2(\fz^-)^{3/2}(x)} \mu\big((\fz^-)^{3/2}\mathbbold{1}_{(-M,\, x)}\big)\!\!<\!0\;\text{ and }\;
\frac{e^{-C(y)}}{2 (\fz^+)^{3/2}(y)}\,
\mu\Big((\fz^+)^{3/2} \mathbbold{1}_{(y,\, N)}\Big)\!\!>\!0.$$
The unique solution $\uz$ to (\ref{08}), or equivalently (\ref{12}), should satisfy
\be{\begin{matrix}
\displaystyle\lim_{x\to -M}\frac{\mu\big((\fz^-)^{3/2} \mathbbold{1}_{(-M,\, x)}\big)}{\sqrt{\fz^-(x)}}\,
\!+\! \mu\big(\fz^-\, \mathbbold{1}_{(-M,\, \uz)}\big)
\!\ge\! \frac{\mu\big((\fz^+)^{3/2} \mathbbold{1}_{(\uz,\, N)}\big)}{\sqrt{\fz^+(\uz)}}\,
\qd\text{and}\\
\displaystyle\lim_{y\to N}\frac{\mu\big((\fz^+)^{3/2} \mathbbold{1}_{(y,\, N)}\big)}{\sqrt{\fz^+(y)}}\,
+ \mu\big(\fz^+\, \mathbbold{1}_{(\uz,\, N)}\big)\ge \frac{\mu\big((\fz^-)^{3/2} \mathbbold{1}_{(-M,\, \uz)}\big)}{\sqrt{\fz^-(\uz)}}.
\end{matrix}}\lb{30}\de

As just mentioned above (cf. \rf{cmf11}{(34)}), we also have
\begin{align}
\max_{z\in [x,\,\uz]} I\!I^-\big(f^{x, y}\big)(z)
&= \max_{z\in [\uz, y]} I\!I^+\big(f^{x, y}\big)(z)
=\frac{h^-(\uz)}{f^{x, y}(x)}=\frac{h^+(\uz)}{f^{x, y}(y)}\nonumber\\
&=\frac{1}{\sqrt{\fz^+(y)}}\,
\mu\Big((\fz^+)^{3/2} \mathbbold{1}_{(y,\, N)}\Big)
+ \mu\big(\fz^+\, \mathbbold{1}_{(\uz,\, y)}\big).\lb{13}
\end{align}
Hence we have arrived at
\begin{align}
&\Big[\sup_{z\in (-M,\, \uz)}I\!I^-(f^{x, y})(z)\Big]
\bigvee \Big[\sup_{z\in (\uz, N)}I\!I^+(f^{x, y})(z)\Big]\nonumber\\
&\qd =\Big[\sup_{z\in (-M,\, x)}I\!I^-(f^{x, y})(z)\Big]
\bigvee \frac{h^+(\uz)}{\sqrt{\fz^+(y)}}
\bigvee \Big[\sup_{z\in (y, N)}I\!I^+(f^{x, y})(z)\Big]\lb{14}
\end{align}
which is also known from \ct{cmf11}. Define
$$I^-(f)(x)=\frac{e^{-C(x)}}{f'(x)}\int_x^{\uz} \frac{e^C}{a}f,\qqd
I^+(f)(x)=-\frac{e^{-C(x)}}{f'(x)}\int_{\uz}^x \frac{e^C}{a}f $$
and
$$\dz_{x,\, \uz}^-\!=\!\!\sup_{z\in (-M,\, x)}\fz_z^-\,\mu (z, \uz),\qqd
\dz_{y,\, \uz}^+\!=\!\!\sup_{z\in (y,\, N)}\fz_z^+\,\mu (\uz, z).
$$
Then we have first by the mean value theorem (both $h^-$ and $f^{x, y}$ are
vanished at $-M$) that
$$\sup_{z\in (-M,\, x)}I\!I^-(f^{x, y})(z)\le \sup_{z\in (-M,\, x)}I^-(f^{x, y})(z)$$
and then by \rf{cmf00}{Lemma 1.2} or \rf{cmf05}{page 97} that
$$\sup_{z\in (-M,\, x)}I^-(f^{x, y})(z)\le 4\,\dz_{x,\, \uz}^-.$$
Here we remark that the supremum in the definition of $\dz_{x,\, \uz}^-$
is taken over $(-M,\, x)$ rather than $(-M,\, \uz)\supset (-M,\, x)$.
Hence the original proof for the last estimate needs a slight modification
using the fact that the function $f^{x, y}$ is a constant on $[x,\,\uz]$.
In parallel, since $h^+$ and $f^{x, y}$ vanish at $N$, we have
$$\sup_{z\in (y,\, N)}I\!I^+(f^{x, y})(z)\le \sup_{z\in (y,\, N)} I^+(f^{x, y})(z)
\le 4\,\dz_{y,\, \uz}^+.$$
Therefore, we have arrived at
\begin{align}
{\underline{\kz}^{\rm DD}}
&\le \inf_{x<y}\bigg\{\Big[\sup_{(-M,\, x)}I\!I^-(f^{x, y})\Big]
\bigvee \frac{h^+(\uz)}{\sqrt{\fz^+(y)}}
\bigvee \Big[\sup_{(y,\, N)}I\!I^+(f^{x, y})\Big]\bigg\}\nonumber\\
&\le \inf_{x<\uz<y}\bigg\{\Big[\sup_{(-M,\, x)}I\!I^-(f^{x, y})\Big]
\bigvee \frac{h^+(\uz)}{\sqrt{\fz^+(y)}}
\bigvee \Big[\sup_{(y,\, N)}I\!I^+(f^{x, y})\Big]\bigg\}\nonumber\\
&\le \inf_{x<\uz<y}\bigg\{\big[4\,\dz_{x,\, \uz}^-\big]\bigvee \frac{h^+(\uz)}{\sqrt{\fz^+(y)}}
\bigvee\big[4\,\dz_{y,\, \uz}^+\big] \bigg\}\nonumber\\
&=: \inf_{x<\uz<y} R(x, y,\,\uz)\nonumber\\
&=:\az .\lb{15}
\end{align}
The restriction $\uz\in (x, y)$ is due to the fact that the eigenfunction of ${{\lz}^{\rm DD}}$
is unimodal and $\uz$ is a mimic of its maximum point.
The use of $I\!I^{\pm}$, $I^{\pm}$ and $\dz^{\pm}$ is now
standard (cf. \ct{cmf05}--\ct{cmf11}, for instance).

We now go to the essential new part of the proof. First,
we claim that for each small $\vz$,
there exist $\bar x\in (-M,\,\uz)$ and $\bar y\in (\uz, N)$ (may depend on $\vz$) such that
\be\fz_{\bar x}^-\,\mu ({\bar x}, \uz)\! \ge\! \frac{R(x_0, y_0, \uz_0)}{4}\!-\!\vz,\qqd
\fz_{\bar y}^+\,\mu (\uz, {\bar y})\!\ge\! \frac{R(x_0, y_0, \uz_0)}{4}\!-\!\vz
\lb{16}\de
for some point $(x_0, y_0, \uz_0)$.
In the present continuous case, the conclusion is clear since the infimum
$\az =R(x^*, y^*, \uz^*)$ is achieved at a point $(x^*, y^*, \uz^*)$ with $x^*\le\uz^*\le y^*$, at which
we have not only $h^-(\uz^*)=h^+(\uz^*)$ but also
\be 4\,\dz_{x^*\!,\, \uz^*}^-=4\,\dz_{y^*\!,\, \uz^*}^+= \frac{h^+(\uz^*)}{\sqrt{\fz^+(y^*)}}.\lb{17}\de
To see this, suppose that at the point $(x, y,\,\uz)$ with $x<\uz<y$, we have
\be \frac{h^+(\uz)}{\sqrt{\fz^+(y)}}> 4 \big[\dz_{x,\,\uz}^- \vee \dz_{y,\,\uz}^+\big].\lb{18}\de
Without loss of generality, assume that $\dz_{x,\,\uz}^- \ge \dz_{y,\,\uz}^+$. We now
fix $y$ and let $\tilde \uz \in (\uz, y]$. Then $\dz_{y,\,\uz}^+\ge \dz_{y, \, \tilde\uz}^+$ by definition. In view of (\ref{13}), we have
$$ \frac{h^+(\uz)}{\sqrt{\fz^+(y)}}>\frac{h^+\big(\tilde\uz\big)}{\sqrt{\fz^+(y)}}.$$
Next, to keep $h^-\big(\tilde\uz\big)=h^+\big(\tilde\uz\big)$, one has a new $\tilde x> x$ by using (\ref{12}) (the left-hand side of (\ref{12}) is decreasing in $x$). Correspondingly, we have $\dz_{\tilde x,\, \tilde\uz}^-\ge \dz_{x,\,\uz}^-$.
In particular, for $\tilde \uz$ closed enough to $\uz$ such that
$$\frac{h^+\big(\tilde\uz\big)}{\sqrt{\fz^+(y)}}\ge 4\,\dz_{\tilde x, \, \tilde\uz}^-,$$
we obtain
$$\aligned
\big[4\,\dz_{x,\,\uz}^-\big] \bigvee \frac{h^+(\uz)}{\sqrt{\fz^+(y)}}
\bigvee \big[4\,\dz_{y,\,\uz}^+\big]
&= \frac{h^+(\uz)}{\sqrt{\fz^+(y)}}\\
&>\frac{h^+\big(\tilde\uz\big)}{\sqrt{\fz^+(y)}}\\
&=\big[4\,\dz_{\tilde x,\, \tilde\uz}^-\big] \bigvee  \frac{h^+\big(\tilde\uz\big)}{\sqrt{\fz^+(y)}}
\bigvee \big[4\,\dz_{y,\, \tilde\uz}^+\big].\endaligned$$
Thus, once (\ref{18}) holds, we can find a new point
$(\tilde x, y, \tilde\uz)$ such that $R(x, y, \uz)> R(\tilde x, y, \tilde\uz)$.
In other words, if the infimum $\az $ is attained at $(x^*, y^*, \uz^*)$, we should have
\be \frac{h^+(\uz^*)}{\sqrt{\fz^+(y^*)}}\le 4 \big[\dz_{x^*\!,\, \uz^*}^- \vee \dz_{y^*\!,\, \uz^*}^+\big].\lb{36}\de
One may handle with the other two cases and finally arrive at (\ref{17}).
Note that instead of (\ref{17}), the following weaker condition is still enough for our purpose.
If at some point $(x, y,\,\uz)$,
\be 4\,\dz_{x,\,\uz}^- =4\,\dz_{y,\,\uz}^+=: \az '\ge  \frac{h^+(\uz)}{\sqrt{\fz^+(y)}},\lb{19}
\de
then we have not only $\az '\ge \az $ but also (\ref{16})
for suitable $\bar x\le x$ and $\bar y\ge y$.
To check (\ref{19}), we first mention that the equation $\dz_{x,\,\uz}^- =\,\dz_{y,\,\uz}^+$ is solvable, at least in the case
that $M, N<\infty$. Because $\dz_{x,\, \uz}^-$ starts
from zero at $x=-M$ and then increases as $x\uparrow$;
$\dz_{y,\,\uz}^+$ also starts
from zero at $y=N$ and then increases as $y\downarrow$.
Therefore, there are a lot of $(x, y)$ satisfying the
required equation.
Next, by (\ref{08}), we can regard $\uz$ as a function of $x$ and $y$.
Then, determine $y$ in terms of $x$ by the equation $\dz_{y,\, \uz(x, y)}^+= \dz_{x,\, \uz(x, y)}^-$.
Now there is only one free variable $x$. We claim that (\ref{19}) holds for some $x$ (and then for some $(x, y,\,\uz)$).
Otherwise, the inverse inequality of (\ref{19}) would hold for all $x$ which contradicts with (\ref{36}).

What we actually need is not the pair $\{x, y\}$ satisfying (\ref{19}) but the pair $\{\bar x, \bar y\}$ satisfying (\ref{16}). From which, the remainder of the proof is very much the same as the one given at the beginning of this section. First, we have
$$(\az /4-\vz)\, \fz^- ({\bar x})^{-1}\le \mu ({\bar x}, \uz),\qqd
(\az /4-\vz)\, \fz^+ ({\bar y})^{-1}\le \mu (\uz, {\bar y}).
$$
Summing up these inequalities, we get
$$(\az /4-\vz)\big[\fz^- ({\bar x})^{-1}+\fz^+ ({\bar y})^{-1}\big]
\le \mu ({\bar x}, {\bar y}).$$
Therefore
$$\aligned
\az /4-\vz &\le \big[\fz^- ({\bar x})^{-1}+\fz^+ ({\bar y})^{-1}\big]^{-1}\mu ({\bar x}, {\bar y})\\
&\le \sup_{x<y}\big[\fz^- ({x})^{-1}+\fz^+ ({y})^{-1}\big]^{-1}\mu ({x}, {y})\\
&=\kz^{\text{\rm DD}}\qd(\text{by } (\ref{05})).
\endaligned$$
Combining this fact with (\ref{15}), we obtain
$$
{\underline{\kz}^{\rm DD}}
\le \az   \nonumber
\le 4 \kz^{\text{\rm DD}}+ 4\,\vz.
$$
Letting $\vz\downarrow 0$,
we have thus proved that ${\underline{\kz}}^{\text{\rm DD}}\le 4 \kz^{\text{\rm DD}}$
as required. The main part of the proof is done since the first Dirichlet
eigenvalue is based on compact sets.

Finally, consider the general case that $M, N\le \infty$.
First, we can rule out the degenerated situation that
$\dz_{x, \uz}^-=\dz_{y, \uz}^+=\infty$.
To see this, rewrite $\kz^{\text{\rm DD}}$ as follows
$$\big(\kz^{\text{\rm DD}}\big)^{-1}= \inf_{x<y}
\big[\big(\fz^- (x)\,\mu(x, y)\big)^{-1}  + \big(\mu (x, y)\, \fz^+ (y)\big)^{-1}\big].$$
It is clear that $\big(\kz^{\text{\rm DD}}\big)^{-1}=0$ and then $\lz^{\text{\rm DD}}=0$ by (\ref{40}).
The corollary becomes trivial. Next, if one of $\dz_{x, \uz}^-$ or $\dz_{y, \uz}^+$
is $\infty$, say $\dz_{x, \uz}^-=\infty$ for instance, then
$$\big(\kz^{\text{\rm DD}}\big)^{-1}=  \Big(\sup_y \,\mu (-M, y) \,\fz^+ (y)\Big)^{-1},$$
i.e.,
$$\kz^{\text{\rm DD}}=  \sup_y\, \mu (-M, y)\, \fz^+ (y).$$
This becomes the essentially known one-side Dirichlet problem.
In the case that both of $\dz_{x, \uz}^-$ and $\dz_{y, \uz}^+ $
are finite,
one may adopt an approximating procedure with finite $M$ and $N$.
This was done in the discrete context, refer to
\rf{cmf10}{Proof of Corollary 7.9 and Proof (c) of Theorem 7.10}.
\deprf

\medskip
To illustrate what was going on in the proof above
and the computation/estimation of ${\underline{\kz}^{\rm DD}}$, we
consider two examples to conclude this section.

\xmp{\bf\rf{cmf11}{Example 5.2}}\qd{\rm Consider the simplest example, i.e. the Laplacian operator on $(0, 1)$.
It was proved in \ct{cmf11} that
$\lz^{\text{\rm DD}}=\pi^2$, $\big(\kz^{\text{\rm DD}}\big)^{-1}=16$ and
$\big({\underline\kz}^{\text{\rm DD}}\big)^{-1}\approx 9.43693$.
The eigenfunction of $\lz^{\text{\rm DD}}$ is $g(x)=\sin(\pi x)$ for which
$g'(1/2)=0$ and so $\bar\uz=1/2$ is the root of equation $\lz_{\uz}^-=\lz_{\uz}^+$.
Because of the symmetry, we have $\uz^*=1/2$ and $y^*=1-x^*$.
Since $\mu=\nu=\d x$, we have $\fz_z^-\, \mu(z, 1/2)=(1/2-z) z$. Thus,
$$\dz_x^-=\sup_{z\in (0, x)}\fz_z^-\, \mu(z, 1/2)=\begin{cases}
(1/2-x)x\qd &\text{ if } x\le 1/4\\
1/16& \text{ if } x\in (1/4, 1/2).
\end{cases}$$
By (\ref{13}) and (\ref{12}), we have
$$\frac{h^-(1/2)}{f^{x,\, 1-x}(x)}=\frac 1 8 -\frac{x^2}{10}.$$
Therefore, each
$$x^*\in \bigg[\frac{20-\sqrt{205}}{78},\; \frac{20+\sqrt{205}}{78}\,\bigg]$$
is a solution to the inequality $4\,\dz_x^-\ge {h^-(1/2)}/{f^{x,\, 1-x}(x)}$.
Correspondingly, we have
$$4\,\dz_{x^*}^-=
\begin{cases}
2 (1-2 x^*)x^*  \qd & \text{if } x^*\in \big[\big(20-\sqrt{205}\,\big)/78,\; 1/4\big]\\
1/4 & \text{if } x^*\in \big[1/4,\; \big(20+\sqrt{205}\,\big)/78\big).
\end{cases}$$
Using this, our conclusion that
$${\underline{\kz}}^{\text{\rm DD}}\le 4\,\dz_{x^*}^-\le 4 \kz^{\text{\rm DD}}$$
can be refined as follows:
$$\big(4\,\kz^{\text{\rm DD}}\big)^{-1}=4\le \big(4\,\dz_{x^*}^-\big)^{-1}
\le \frac{35-\sqrt{41/5}}{4}\approx 8.034< 9.43693\approx\big({\underline\kz}^{\text{\rm DD}}\big)^{-1}.$$
It follows that there are many solutions $x^*$, and so we have a lot of freedom in choosing $(\uz^*, x^*, y^*)$ for (\ref{19}).
However, the maximum of $\big(4\,\dz_{x^*}^-\big)^{-1}$ is attained only at the point $x^*$ which is the
smaller root of equation: $4\,\dz_x^-= {h^-(1/2)}/{f^{x,\, 1-x}(x)}$.
}
\dexmp

The next example is unusual since for which the lower bound $\big(4\,\kz^{\text{\rm DD}}\big)^{-1}$
is sharp. Hence, there is no room for the improvement $\big({\underline\kz}^{\text{\rm DD}}\big)^{-1}$.
The proof above seems rather dangerous for this example since at each step
$$\aligned
\big(\lz^{\text{\rm DD}}\big)^{-1}&\le {\underline\kz}^{\text{\rm DD}}
\le \text{Est}(I^{\pm}(f))\le \text{Est}(\dz_{\cdot,\, \uz}^{\pm})\\
&\le 4\big(\ffz (\bar x, \bar y)+\vz\big)
\le 4 \sup_{x<y}\big(\ffz (x, y)+\vz\big)=4\big( {\kz}^{\text{\rm DD}}+\vz\big)
\endaligned$$
for some $\ffz$, where Est\,$(H)$ means the estimate using $H$,
one may lose something. Here we have also explained the reason
why ${\underline{\kz}^{\rm DD}}$ is often much better than
$4\,{\kz}^{\rm DD}$ as shown in the last example.

\xmp{\bf\rf{cmf11}{Example 5.3}}\qd{\rm Consider the operator
$L=\d^2/\d x^2 +b \d/\d x$ with $b>0$ on $(0, \infty)$. It was checked in \ct{cmf11} that
$\lz^{\rm DD}=b^2/4$, $\big(\kz^{\rm DD}\big)^{-1}=b^2$ and
so the lower estimate $\big(4\,\kz^{\rm DD}\big)^{-1}$ is sharp.
The eigenfunction of $\lz^{\text{\rm DD}}$ is $g(x)=x e^{-b x/2}$ for which
$g'(2/b)=0$ and so $\bar\uz=2/b$ solves the equation $\lz_{\uz}^-=\lz_{\uz}^+$.
We have $C(x)= bx$, $\mu (\d x)=e^{b x} \d x$,
$$\fz^-(s)=\int_0^s e^{-b z}\d z=\frac{1}{b}\big(1- e^{-bs}\big)\qd\text{and}\qd
\fz^+(s)=\int_s^\infty e^{-bz}\d z= \frac{1}{b}e^{-bs}.$$
We begin our study on the equation $\dz_{x,\,\uz}^-= \dz_{y,\,\uz}^+$ rather than
Eq.(\ref{12}) since the former one is simpler. Note that the function
$$\fz_z^-\,\mu (z, \uz)=\frac 1 {b^2} \big(1 - e^{-b z}\big) \big(e^{b \uz}-e^{b z}\big),\qqd z\in (0, \uz]$$
achieves its maximum $b^{-2}\big(e^{b\uz/2}-1\big)^2$
at $z=\uz/2$ and the function
$$\mu (\uz, z)\fz_z^+=\frac 1 {b^2} \big(1-e^{b(\uz - z)} \big),\qqd z\ge \uz $$
achieves its maximum $1/b^2$ at $\infty$. Hence
$$\dz_{x,\,\uz}^-= \frac{1}{b^2}\big(e^{b\uz/2}-1\big)^2\qd \forall x\in [\uz/2, \uz]\qd\text{and}\qd
\dz_{y,\,\uz}^+=\frac 1 {b^2}\qd\forall y\ge \uz.$$
Solving the equation
$$\frac{1}{b^2}\big(e^{b\uz/2}-1\big)^2=\frac 1 {b^2},$$
we get
$\uz^*= 2 b^{-1} \log 2.$
To study (\ref{30}), note that
$$\aligned
&\frac{1}{\sqrt{\fz^+(y)}}
\mu\Big((\fz^+)^{3/2} \mathbbold{1}_{(y,\, \infty)}\Big)
+\mu\big(\fz^+ \mathbbold{1}_{(\uz,\, y)}\big)=\frac{2}{b^2}+\frac 1 b (y-\uz)\\
&\frac{1}{\sqrt{\fz^-(x)}}\mu\Big(\!(\fz^-)^{3/2} \mathbbold{1}_{(0,\, x)}\!\Big)
\!+\!\mu\big(\fz^- \mathbbold{1}_{(x,\,\uz)}\big)\\
&\qqd\qqd= \frac{1}{b^2}\bigg\{2-b\uz+e^{b\uz}+bx
-\frac{3\big(b x +\log\big(1+\sqrt{1-e^{-bx}}\,\big)\big)}{2 \sqrt{1-e^{-bx}}}
\bigg\}.
\endaligned$$
Then the second inequality in (\ref{30}) is trivial and the first one there becomes
$$
\frac{1}{b^2}\big(2-b\uz+e^{b \uz}\big)\ge \frac{2}{b^2}.$$
It is now easy to check that $\uz^*=2 b^{-1} \log 2$ does not satisfy this inequality.
In other words, there is no required solution $(x^*, y^*, \uz^*)$ under the restriction $x^*\in [\uz^*/2, \uz^*]$.
Thus, unlike the last example, there is not much freedom in choosing $(x^*, y^*, \uz^*)$ for (\ref{19}).
However, this does not finish the story since the solution $x^*$ may belong to $[0, \uz^*/2)$.

We are now looking for a solution $x^*$ in the interval $[0, \uz^*/2)$.
When $x\le \uz/2$, the maximum of the function
$\sup_{z\le x} \fz_z^-\,\mu (z, \uz)$
on $[0, x]$ is achieved at $x$. Hence
$$\dz_{x,\,\uz}^-= \frac 1 {b^2} \big(1 - e^{-b x}\big) \big(e^{b \uz}-e^{b x}\big)\qd \forall x\in (0, \, \uz/2]\qd\text{and}\qd
\dz_{y,\,\uz}^+=\frac 1 {b^2}\qd\forall y\ge \uz.$$
Solving the equation
$$\frac 1 {b^2} \big(1 - e^{-b x}\big) \big(e^{b \uz}-e^{b x}\big)=\frac 1 {b^2},$$
we obtain
$\uz^*=x-b^{-1} \log\big(1-e^{-b x}\big).$
Besides, solving the equation $4\,\dz_{y^*\!, \, \uz^*}^+=h^+(\uz^*)\big/\!\sqrt{\fz^+(y^*)}\,$,
we get $y^* = 2/b + \uz^*$. Inserting these into Eq.(\ref{12}), we obtain
$$\frac{e^{2 b x}}{e^{b x}-1}+\log \left(1-e^{-b x}\right)=2+
\frac{3}{2 \sqrt{1-e^{-b x}}} \left(b x+2 \log \left(\sqrt{1-e^{-b x}}+1\right)\right).$$
From this, we obtain the required solution $x^*$ as shown by Figures 1 and 2 below,
noting that the constraint that
$x^*\le \uz^*/2$ is equivalent to $x^*\le b^{-1}\log 2$. Having $x^*$ at hand, it is clear that
the solution $\uz^*$ here is very different from $2/b$.
\begin{center}{{{\includegraphics{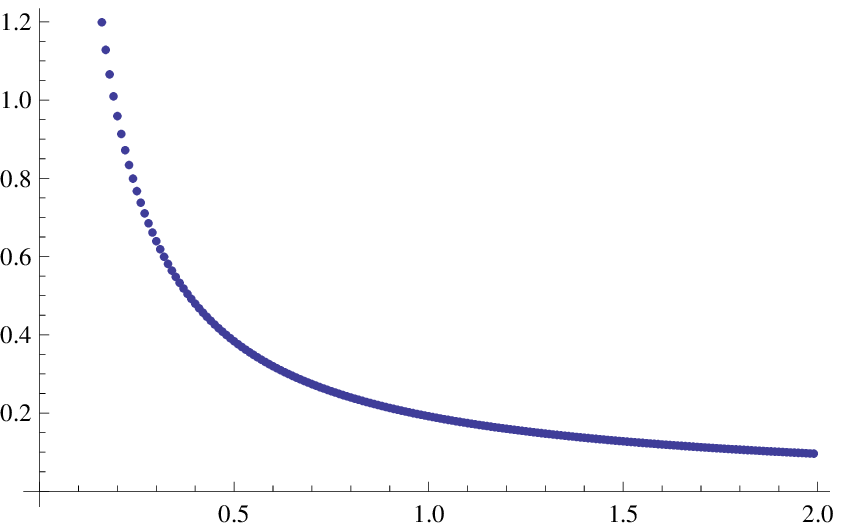}\hskip-4.8truecm}
}
\vskip-4.2truecm{\hskip-1.8truecm\includegraphics{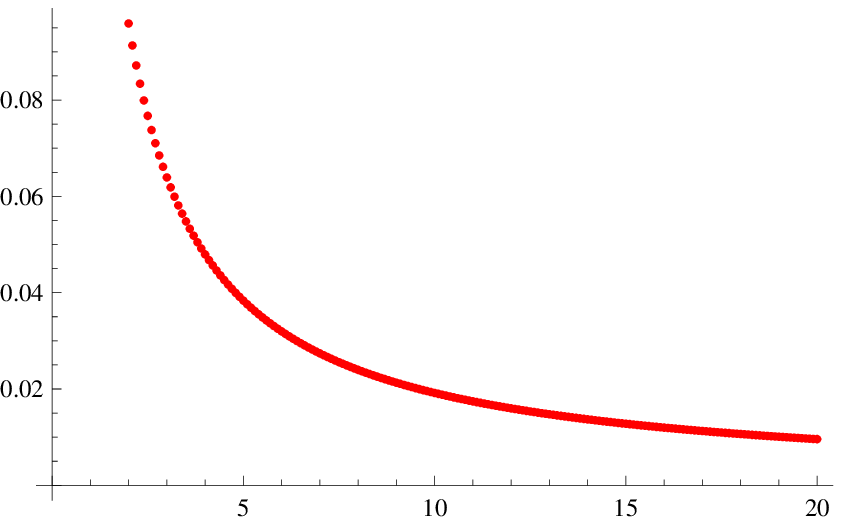}}\newline
{\bf Figure 1--2}\qd\rm Solution of $x^*=x^*(b)$ when $b$ varies on $(0, 2]$
(the curve on right) and on $[2, 20]$ (the curve on left), respectively.}\end{center}

To see that the solutions $(\bar x, \bar y)$ to (\ref{16}) may not be unique, keeping $\uz^*$ to be the same as in the last paragraph
but replace $y^*$ with a smaller one $\bar y = b^{-1}+\uz^*$, then one can find a point $\bar x$
satisfying Eq.\,(\ref{12}).
}
\dexmp

\section{Birth--death processes (discrete case)}

This section deals with the discrete case which is parallel in principal to the
continuous one studied above, but it is quite involved and so is worth to write down some details here.

The state space is
$$E=\{i\in {\mathbb Z}: -M-1< i< N+1\},\qqd M, N\le \infty.$$ The transition rates
$Q=(q_{ij})$ are as follows: $b_i:=q_{i,i+1}>0$, $a_i:=q_{i, i-1}>0$, $q_{ii}=-(a_i+b_i)$, $i\in E$.
$q_{ij}=0$ for other $i\ne j$.
Thus, we have $a_{-M}>0$ if $M<\infty$ and similarly for $b_N$. The operator of the process becomes
$$\ooz f(i)=b_i\big(f_{i+1}-f_i\big)+ a_i\big(f_{i-1}-f_i\big), \qqd i\in E$$
with a convention $f_{-M-1}=0$ if $M<\infty$ and  $f_{N+1}=0$ if $N<\infty$.
Next, define the speed (or invariant, or symmetric) measure $\mu$ as follows.
Fix a reference point $o\in E$ and set
\begin{gather}
\mu_{o+n}=\frac{a_{o -1} a_{o -2}\cdots a_{o +n+1}}
 {b_{o } b_{o-1 }\cdots b_{o +n}}, \qqd -M-1-o< n\le -2,\nonumber\\
\mu_{o-1 }=\frac{1}{b_{o } b_{o -1}},\qqd
\mu_{o }=\frac{1}{a_{o } b_{o }},\qqd
 \mu_{o+1 }=\frac{1}{a_{o } a_{o +1}},\nonumber\\
 \mu_{o+n }=\frac{b_{o +1} b_{o +2}\cdots b_{o +n-1}}
 {a_{o } a_{o+1 }\cdots a_{o +n}}, \qquad  2\le n<N+1-o.\nonumber
\end{gather}
A change of the reference point $o$ leads to a constant factor only to the sequence
$(\mu_i)$ and so does not make any influence to the results below.
Corresponding to $\ooz$, the Dirichlet form is
\begin{gather}
D(f)=\sum_{-M-1<i\le o} \mu_i a_i (f_i-f_{i-1})^2
        +\sum_{o\le i <N+1} \mu_i b_i (f_{i+1}-f_{i})^2,\nonumber\\
  f\in {\scr K},\; f_{-M-1}=0\;\text{if }M<\infty \text{ and }f_{N+1}=0\;\text{if }N<\infty, \nonumber
\end{gather}
where ${\scr K}$ is the set of functions on $E$ with compact supports. Having these preparations
at hand, one can define the eigenvalues $\lz^{\text{\rm DD}}$ and  $\lz^{\text{\rm NN}}$ on $L^2(\mu)$
as in the first section.

To state our main result in this context, we need more notation. Define
$${\scr C}_+=\{f|_E>0: f_{-M-1}=0\;\text{if }M<\infty \text{ and }f_{N+1}=0\;\text{if }N<\infty\}.$$
Given $f\in {\scr C}_+$, define $h^{\pm}=h_f^{\pm}$ as follows.
$$\aligned
h_i^-&=\sum_{k=-M}^i \frac{1}{\mu_k a_k} \sum_{\ell=k}^{\uz} \mu_{\ell} f_{\ell}
=\sum_{\ell=-M}^{\uz} \mu_{\ell} f_{\ell} \fz_{\ell\wedge i}^-, \qd i\le \uz,\\
h_i^+&=\sum_{k=i}^N \frac{1}{\mu_k b_k} \sum_{\ell=\uz}^{k} \mu_{\ell} f_{\ell}
=\sum_{\ell=\uz}^{N} \mu_{\ell} f_{\ell} \fz_{\ell\vee i}^+, \qd i\ge \uz,
\endaligned$$
where
$$\fz_i^-=\sum_{k=-M}^i \frac{1}{\mu_k a_k},\qqd \fz_k^+=\sum_{\ell=k}^N \frac{1}{\mu_{\ell} b_{\ell}}$$
and $\uz\in (-M-1, N+1)$ will be specified soon.
Applying $h_f^{\pm}$ to the test function $f=f^{m, n}\,(m, n\in E,\, m\le n)$:
$$f_i^{m,n}=
\begin{cases} \sqrt{\fz_n^+ \fz_{i\wedge m}^-/\fz_m^-}\qqd & i\le n\\
   \sqrt{\fz_i^+} & i\ge n,
\end{cases}$$
we obtain a condition for $\uz$ which is an analog of (\ref{30}):
\be{\begin{matrix}
\displaystyle\lim_{m\to -M}\frac{1}{\sqrt{\fz_m^-}}\!
\sum_{k=-M}^{m-1}\! (\fz_k^-)^{3/2}\mu_k \!+\!\! \sum_{k=-M}^{\uz}\! \fz_k^-\mu_k
\displaystyle\!\ge\! \frac{1}{\sqrt{\fz_{\uz}^+}}\!
\sum_{k=\uz}^N  (\fz_k^+)^{3/2}\mu_k\qd \text{and}\\
\displaystyle\lim_{n\to N}\frac{1}{\sqrt{\fz_n^+}}\!
\sum_{k=n+1}^N\!  (\fz_k^+)^{3/2}\mu_k + \sum_{k=\uz}^N  \fz_k^+\mu_k
\displaystyle\!\ge\! \frac{1}{\sqrt{\fz_{\uz}^-}}\!
\sum_{k=-M}^{\uz} (\fz_k^-)^{3/2}\mu_k.
\end{matrix}}\lb{32}\de
However, in the discrete situation, one can not expect (\ref{08}). This leads to a
serious change. To explain the main idea, let us return to Theorem \ref{t21}.
Because the derivative of the eigenfunction of $\lz^{\rm DD}$ has uniquely one zero point,
say $\uz$. We can split the interval $(-M, N)$ into two parts having a common boundary
$\uz$. Thus, the original process is divided into processes having a common reflecting
boundary $\uz$. Theorem \ref{t21} says that the original $\lz^{\rm DD}$ can be represented
by using the principal eigenvalues of these sub-processes. This idea is the starting point of \ct{czz03}, as already used in the first proof in Section 2. Since the maximum point $\uz$ is unknown in advance, in the
original formulation, $\uz$ is free and then there is an additional term $\sup_{\uz}$
in the expression of Theorem \ref{t21}. This term was removed in \ct{cmf11}, choosing
$\uz$ as a mimic of the maximum point of the eigenfunction. Unfortunately, such a mimic still does not work in the
discrete case, we may lost (\ref{08}) and more seriously,
the eigenfunction may be a simple echelon but not a unimodal (cf. \rf{cmf10}{Definition 7.13}).
Therefore, more work is required. Again, the idea goes back to \ct{czz03} except
here the choice of $\uz$ is based on (\ref{32}).
The first key step of the method is constructing two birth--death processes on the left- and
the right-hand sides, separately. As before, the two processes have Dirichlet boundaries
at $-M-1$ and $N+1$ but they now have a common Neumann boundary at $\uz\in E$.
Let us start from the birth--death process with rates $(a_i, b_i)$ and state space $E$.
Fix a constant $\gz>1$.
\begin{itemize}
\item [(L)] The process on the left-hand side has state space $E^{\uz -}=\{i: -M-1<i\le \uz\}$,
reflects at $\uz$ (and so $b_{\uz}=0$).
Its transition structure is the same as the original one except $a_{\uz}$ is replaced by $a_{\uz}^{-, \gz}:=\gz a_{\uz}$.
Then for this process, the sequence $(\mu_i: i\in E^{\uz -})$ is the same as the original one except the original
$\mu_{\uz}$ is replaced by $\mu_{\uz}/\gz$. Hence, the sequence $(\mu_i a_i: i\in E^{\uz -})$
keeps the same as original.
\item [(R)] The process on the right-hand side has state space $E^{\uz +}=\{i: \uz\le i <N+1\}$,
reflects at $\uz$ (and then $a_{\uz}=0$). Its transition structure is again
the same as the original one except $b_{\uz}$ is replaced by $b_{\uz}^{+, \gz}:=\gz (\gz-1)^{-1} b_{\uz}$.
Then for this process, the sequence $(\mu_i: i\in E^{\uz +})$ is the same as the original one except the original
$\mu_{\uz}$ is replaced by $(1-\gz^{-1})\mu_{\uz}$. Hence, the sequence $(\mu_i b_i: i\in E^{\uz +})$
remains the same as original.
\end{itemize}
Noting that $a_{\uz}^{-, \gz} \downarrow a_{\uz}$ and $b_{\uz}^{+, \gz}\uparrow \infty$ as $\gz\downarrow 1$,
$a_{\uz}^{-, \gz} \uparrow \infty$ and $b_{\uz}^{+, \gz}\downarrow b_{\uz}$ as $\gz\uparrow \infty$,
the constant $\gz$ plays a balance role for the principal eigenvalues of these processes.
From here, following the first proof given in Section 2 and using \ct{czz03} and \rf{cmf10}{Theorem 7.10}, one
can prove the basic estimate
$\lz^{\text{\rm DD}}\ge \big(4\kz^{\text{\rm DD}}\big)^{-1}$
in the present context. Certainly, the parallel proof works also in the ergodic case.

We now continue our study on the discrete analog of Corollary \ref{t11}\,(1).
The quantity $\fz^{\pm}$ needs no change. But $h^{\pm}$ has to be modified as follows.
$$\aligned
h_i^{-, \gz}&=\sum_{k=-M}^i \frac{1}{\mu_k a_k} \bigg[\sum_{k\le \ell\le\uz-1} \mu_{\ell} f_{\ell}+\frac{1}{\gz}\mu_{\uz}f_{\uz}\bigg], \qqd i\le \uz,\\
h_i^{+, \gz}&=\sum_{k=i}^N \frac{1}{\mu_k b_k} \bigg[\frac{\gz-1}{\gz}\mu_{\uz}f_{\uz}+
\sum_{\uz+1\le\ell\le k} \mu_{\ell} f_{\ell}\bigg], \qqd i\ge \uz.
\endaligned$$

Finally, define $I\!I^{\pm, \gz}(f)=h^{\pm, \gz}/f$.
It is now more convenient to write the test functions on $E_{\uz}^{\pm}$ separately:
$$
f_i^{-, m}=\sqrt{\fz_{i\wedge m}^-},\qd i\le \uz,\qqd
f_i^{+, n}=\sqrt{\fz_{i\vee n}^+},\qd i\ge \uz.
$$
Comparing with the original $f^{m, n}$, here a factor acting on $f^{-, m}$ is ignored
\big(the reason why one needs the factor in the original case is for $f_{\uz}^{-, m}=f_{\uz}^{+, n}$\big).

\crl\lb{t31} {\cms We have
\be \lz^{\rm DD}\ge \big({\underline\kz}^{\text{\rm DD}}\big)^{-1}\ge \big(4\,\kz^{\text{\rm DD}}\big)^{-1},\lb{20}\de
where
\begin{gather}
\mbox{\hspace{-1em}}{\underline\kz}^{\text{\rm DD}}\!=\inf_{\uz: \text{\rm(\ref{32}) holds}}\,\inf_{m\le \uz\le n}\inf_{\gz>1}\!\Big\{\!\Big[\sup_{ E\owns i\le\uz} I\!I_i^{-, \gz}(f^{-, m})\Big]
\!\bigvee\! \Big[\sup_{\uz\le i\in E} I\!I_i^{+, \gz}(f^{+, n})\Big]\!\Big\},\\
\mbox{\hspace{-2em}}\big(\kz^{\text{\rm DD}}\big)^{-1}\!=\inf_{m, n\in E:\; m \le n}\bigg[\bigg(\sum_{i=-M}^m\frac{1}{\mu_i a_i}\bigg)^{-1}\!
\!\!+ \bigg(\sum_{i=n}^{N}\frac{1}{\mu_i
b_i}\bigg)^{-1}\bigg]\bigg(\sum_{j=m}^n \mu_j\bigg)^{-1}\!\!\!.
\end{gather}
}
\decrl

\prf By using an approximating procedure, one may assume that $M, N<\infty$
(cf. \rf{cmf10}{Proof of Corollary 7.9 and Proof (c) of Theorem 7.10}).
Fix $\uz\in [m, n]$ and define
\begin{gather}
I_i^{-, \gz}(f)=\frac{1}{\mu_i a_i (f_i-f_{i-1})}\bigg[\frac{1}{\gz}\mu_{\uz}+\sum_{i\le \ell\le \uz-1} \mu_{\ell}\bigg],\qqd i\le \uz\nonumber\\
I_i^{+, \gz}(f)=\frac{1}{\mu_i b_i (f_i-f_{i+1})}\bigg[\frac{\gz-1}{\gz}\mu_{\uz}+\sum_{\uz+1\le \ell \le i} \mu_{\ell}\bigg],\qqd i\ge \uz\nonumber\\
\dz_{m,\, \uz}^{-, \gz}=\sup_{i\le m}\fz_i^- \bigg[\frac{1}{\gz}\mu_{\uz}+\!\!\sum_{i\le \ell\le \uz-1} \mu_{\ell}\bigg], \qqd
\dz_{n,\, \uz}^{+, \gz}=\sup_{i\ge n} \fz_i^+\bigg[\frac{\gz-1}{\gz}\mu_{\uz}+\!\!\sum_{\uz+1\le \ell \le i} \mu_{\ell}\bigg].\nonumber
\end{gather}
We have
\begin{align}
\frac{h_{\uz}^{-, \gz}}{f_m^{-, m}}
&=\frac{1}{\sqrt{\fz_m^-}}\sum_{k=-M}^{m-1}  (\fz_k^-)^{3/2}\mu_k
+ \sum_{k=m}^{\uz-1}  \fz_k^- \mu_k +\frac{1}{\gz} \fz_{\uz}^- \mu_{\uz},\\
\frac{h_{\uz}^{+, \gz}}{f_n^{+, n}}
&=\frac{1}{\sqrt{\fz_n^+}}\sum_{k=n+1}^N (\fz_k^+)^{3/2}\mu_k
+ \sum_{k=\uz+1}^{n} \fz_k^+\mu_k +\frac{\gz-1}{\gz} \fz_{\uz}^+ \mu_{\uz}.\end{align}
For simplicity, let
$$H(m, n,\,\uz, \gz)= \max\bigg\{\frac{h_{\uz}^{-, \gz}}{f_m^{-, m}}\mathbbold{1}_{\{m<\uz\}},\;
\frac{h_{\uz}^{+, \gz}}{f_n^{+, n}}\mathbbold{1}_{\{\uz<n\}}\bigg\}.$$
By \rf{cmf10}{Theorem 7.10\,(1), Sections 4, 2 and 3}, we obtain
\begin{align}
\big({{\underline\kz}^{\rm DD}}\big)^{-1}\!\!
&\le\! \inf_{\begin{subarray}{c}\gz>1 \\ [m, n]\owns \uz\\ \uz: \text{\rm(\ref{32}) holds}\end{subarray}}\!\bigg\{\!\bigg[\bigvee_{i\le m}\! I\!I_i^{-, \gz}(f^{-, m})\bigg]
\!\bigvee\! H(m, n,\,\uz, \gz)
\!\bigvee\! \bigg[\bigvee_{i\ge n}\! I\!I_i^{+, \gz}(f^{+, n})\bigg]\!\bigg\}\nonumber\\
&\le \inf_{\begin{subarray}{c}\gz>1 \\ [m, n]\owns \uz\\ \uz: \text{\rm(\ref{32}) holds}\end{subarray}}\bigg\{\!\bigg[\bigvee_{i\le m}\! I_i^{-, \gz}(f^{-, m})\bigg]
\!\bigvee\!H(m, n,\,\uz, \gz)
\!\bigvee\! \bigg[\bigvee_{i\ge n}\! I_i^{+, \gz}(f^{+, n})\bigg]\!\bigg\}\nonumber\\
&\le \inf_{\begin{subarray}{c}\gz>1 \\ [m, n]\owns \uz\\ \uz: \text{\rm(\ref{32}) holds}\end{subarray}}\bigg\{\Big[ 4 \dz_{m, \,\uz}^{-, \gz}\Big]
\bigvee H(m, n,\,\uz, \gz)
\bigvee \Big[4 \dz_{n,\, \uz}^{+, \gz}\Big]\bigg\}\nonumber\\
&=: \inf_{\uz: \text{\rm(\ref{32}) holds}}\,\inf_{[m, n]\owns \uz}\,\inf_{\gz>1}\, R(m, n,\,\uz, \gz)\nonumber\\
&=:\az. \lb{24}
\end{align}
The point we need two terms in the expression of $H$,
rather than one only in the continuous case, is the loss of an analog of (\ref{08}): here we may not have
$h_{\uz}^{-, \gz}=h_{\uz}^{+, \gz}$.
We now choose a candidate of $\uz^*$ (independent of $m$, $n$) from (\ref{32}) and then choose $\{m^*, n^*\}$ with
$m^*\le \uz^*\le n^*$ (may not be unique) so that $(m^*, n^*, \uz^*)$ satisfies the following inequalities
\be{\begin{matrix}\displaystyle
\Big(\dz_{m,\, \uz}^{-, 1}=\!\Big)\;\sup_{i\le m}\fz_i^- \bigg[\mu_{\uz}+\!\!\sum_{i\le \ell\le \uz-1} \mu_{\ell}\bigg]\ge
\sup_{i\ge n} \fz_i^+ \sum_{\uz+1\le \ell \le i} \mu_{\ell}\;\Big(\!= \dz_{n,\, \uz}^{+, 1}\Big)\qd\text{and}\\
\displaystyle\Big(\dz_{m,\, \uz}^{-, \infty}=\!\Big)\;\sup_{i\le m}\fz_i^- \sum_{i\le \ell\le \uz-1} \mu_{\ell} \le
\sup_{i\ge n} \fz_i^+\bigg[\mu_{\uz}+\!\!\sum_{\uz+1\le \ell \le i} \mu_{\ell}\bigg] \;\Big(\!=\dz_{n,\, \uz}^{+, \infty}\Big).\end{matrix}}
\lb{38}\de
Roughly speaking, the condition (\ref{08}) in the continuous case is replaced by a much weaker one (\ref{32}) and the condition $\dz_{m, \,\uz}^{-}= \dz_{n,\, \uz}^{+}$ is replaced by (\ref{38}).
Instead, let $\gz^*$ be the unique solution to the equation
$$\dz_{m^*, \,\uz^*}^{-, \gz}= \dz_{n^*,\, \uz^*}^{+, \gz},\qqd \gz\in [1, \infty]$$
for each fixed pair $\{m^*, n^*\}$: $m^*\le\uz^*\le n^*$.
Here is the balance role played by $\gz$ as mentioned before. As an analog of the continuous case, we are interested in those $\{m^*, n^*\}$: $[m^*\!,\,n^*]\owns \uz^*$ having the property
\be\dz_{m^*\!, \,\uz^*}^{-, \gz^*}= \dz_{n^*\!,\, \uz^*}^{+, \gz^*}\ge\frac{1}{4} \max\bigg\{\frac{h_{\uz^*}^{-, \gz^*}}{f_{m^*}^{-, m^*}}\mathbbold{1}_{\{m^*<\uz^*\}},\;
\frac{h_{\uz^*}^{+, \gz^*}}{f_{n^*}^{+, n^*}}\mathbbold{1}_{\{\uz^*<n^*\}}\bigg\}.\lb{34}\de
Unlike the continuous case, here we may have to repeat the procedure in choosing $(m^*, n^*, \gz^*)$ since $\uz^*$ suggested by (\ref{32}) may not be unique.
Note that the right-hand side of (\ref{34}) is trivial in the particular case that $m^*=n^*=\uz^*$.
Thus, for sufficiently small $\vz>0$, we may choose $(\bar m, \bar n)$ with $[\bar m, \bar n]\owns \uz^*$ and
${\bar\gz}\in (1, \infty)$ such that
\begin{align}
&\fz_{\bar m}^-\bigg[\frac{1}{{\bar\gz}}\mu_{\uz^*}+ \sum_{\ell=\bar m}^{\uz^*-1} \mu_{\ell}\bigg]\ge \frac{R(m^*, n^*, \uz^*, \gz^*)}{4}-\vz \qd\text{and}\nonumber\\
  &\fz_{\bar n}^+\bigg[\frac{{\bar\gz}-1}{{\bar\gz}}\mu_{\uz^*}+\sum_{\ell=\uz^*+1}^{\bar n} \mu_{\ell} \bigg]\ge \frac{R(m^*, n^*, \uz^*, \gz^*)}{4}-\vz. \nonumber\end{align}
Therefore, we have
$$\bigg(\frac{\az}{4}-\vz\bigg)\big(\fz_{\bar m}^-\big)^{-1}\!\!\le
   \frac{1}{{\bar\gz}}\,\mu_{\uz^*}+\! \sum_{\ell=\bar m}^{\uz^*-1} \mu_{\ell},\qqd
\bigg(\frac{\az}{4}-\vz\bigg)\big(\fz_{\bar n}^+\big)^{-1}\!\!\le
   \frac{{\bar\gz}-1}{{\bar\gz}}\,\mu_{\uz^*}+\!\sum_{\ell=\uz^*+1}^{\bar n}\! \mu_{\ell}.
$$
Summing up these inequalities, it follows that
$$\aligned
\bigg(\frac{\az}{4}\!-\!\vz\bigg)
\Big\{\big(\fz_{\bar m}^-\big)^{-1}\!\!+\!
\big(\fz_{\bar n}^+\big)^{-1}\Big\}
\!\le\! \sum_{\ell=\bar m}^{\uz^*-1}\! \mu_{\ell}+\frac{1}{{\bar\gz}}\mu_{\uz^*}\!+\!
\frac{{\bar\gz}-1}{{\bar\gz}}\mu_{\uz^*}\!+\!\!\sum_{\ell=\uz^*+1}^{\bar n}\! \mu_{\ell}
\! =\!\sum_{\ell=\bar m}^{\bar n} \mu_{\ell}.
\endaligned$$
The remainder of the proof is the same as in the continuous situation.\deprf

The following example is almost the simplest one but is indeed very helpful to understand Corollary \ref{t31}
and its proof.

\xmp{\bf \rf{czz03}{Example 2.3} and \rf{cmf10}{Example 7.6\,(2)}}\qd{\rm Let $M=-1$, $N=2$, $b_1=1$, $b_2=2$,
$$a_1=\frac{2-\vz^2}{1+\vz}, \qqd \vz\in \big[0, \sqrt{2}\,\big) \qqd\text{and\qqd $a_2=1$}.$$
Then ${{\lz}^{\rm DD}}=2-\vz$.
It is known that
$$\kz^{\rm DD}=\frac{1}{\lz_0}
-\begin{cases}
{\vz^2}(8 - 4\, \vz^2 + \vz^3)^{-1}\qd &\text{if } \vz\in \big[0,\; \big(\sqrt{13}-1\big)/{3}\big]\\
(8 + 2\, \vz - 3\, \vz^2)^{-1} &\text{if } \vz\in \big[\big(\sqrt{13}-1\big)/{3},\; \sqrt{2}\,\big).
\end{cases}$$

We are now going to compare ${\underline\kz}^{\rm DD}$ with $4\,\kz^{\rm DD}$.
First, we have $\mu_1=\mu_2=1$, $\mu_1 a_1=a_1$, $\mu_1 b_1=b_1$ and $\mu_2 b_2=b_2$.
Next, (\ref{32}) holds for a small part of $\vz$ when $\uz=1$ but never holds if $\uz=2$. Hence, we choose
$\uz=1$. Then $m=1$ and
$$\fz_1^-= \frac{1}{a_1}=\frac{1+\vz}{2-\vz^2},\qqd
\fz_1^+=\frac{1}{b_1}+\frac{1}{b_2}=\frac 3 2, \qqd \fz_2^+=\frac{1}{b_2}=\frac 1 2.$$
Furthermore
$$\dz_{m,\,\uz}^{-, \gz}=\frac{1}{\gz a_1}=\frac{1+\vz}{\gz(2-\vz^2)}.$$
By (\ref{38}), we have $n=1$ or $2$.

(1) When $n=1$, we have
$$\dz_{n,\,\uz}^{+, \gz}
=\bigg[\fz_1^+\bigg(1-\frac{1}{\gz}\bigg)\bigg]\bigvee
\bigg[\fz_2^+\bigg(2-\frac{1}{\gz}\bigg)\bigg]
=\begin{cases}
\displaystyle\frac{3}{2}\bigg(1-\frac{1}{\gz}\bigg)\qd &\text{if } \gz\ge 2\\
\displaystyle 1-\frac{1}{2\gz} &\text{if } \gz\in (1, 2).
\end{cases}$$
Clearly, the equation $\dz_{m,\,\uz}^{-, \gz}=\dz_{n,\,\uz}^{+, \gz}$ has a unique solution
$$\gz=\begin{cases}
\displaystyle\frac{8 + 2 \vz - 3 \vz^2}{3 (2 - \vz^2)}\ge 2 \qd &\text{if $\displaystyle \vz\in\bigg[\frac{\sqrt{13}-1}{3},\, \sqrt{2}\,\bigg)$}\\
\displaystyle\frac{4 + 2 \vz - \vz^2}{2 (2 - \vz^2)}\in (1, 2) \qd &\text{if $\displaystyle \vz\in\bigg(0,\, \frac{\sqrt{13}-1}{3}\bigg)$}.
\end{cases}$$
Correspondingly, with $m=n=\uz=1$, we have
$$\dz_{m,\,\uz}^{-, \gz}=\dz_{n,\,\uz}^{+, \gz}=
\begin{cases}
\displaystyle\frac{3(1+\vz)}{8 + 2 \vz - 3 \vz^2}\qd &\text{if $\displaystyle \vz\in\bigg[\frac{\sqrt{13}-1}{3},\, \sqrt{2}\,\bigg)$}\\
\displaystyle\frac{2(1+\vz)}{4 + 2 \vz - \vz^2}\qd &\text{if $\displaystyle \vz\in\bigg(0,\, \frac{\sqrt{13}-1}{3}\bigg)$}.
\end{cases}$$
It is interesting that the last quantity coincides with $4\,\kz^{\rm DD}$. We have thus arrived at (\ref{34}) since we are in the particular case: $m^*=n^*=\uz^*$.

(2) When $n=2$, we have
$$
\frac{h_{\uz}^{+, \gz}}{f_n^{+, n}}
=\fz_2^+ +\frac{\gz-1}{\gz} \fz_1^+=2-\frac{3}{2\gz},\qqd
\dz_{n,\,\uz}^{+, \gz}
=1-\frac{1}{2\gz}.$$
Clearly,
$$\frac{h_{\uz}^{+, \gz}}{f_n^{+, n}}\le 4\,\dz_{n,\,\uz}^{+, \gz}\qd\text{iff $\gz\ge 1/4$}.$$
As we have seen above, the solution to the equation $\dz_{m,\,\uz}^{-, \gz}=\dz_{n,\,\uz}^{+, \gz}$ is
$$\gz=\frac{4 + 2 \vz - \vz^2}{2 (2 - \vz^2)}>1 \qd \text{on $\big(0, \sqrt{2}\,\big)$}.$$
Then
$$\dz_{m,\,\uz}^{-, \gz}=\dz_{n,\,\uz}^{+, \gz}=\frac{2(1+\vz)}{4 + 2 \vz - \vz^2}>\frac{h_{\uz}^{+, \gz}}{4 f_n^{+, n}}\qd
(m=\uz=1,\; n=2).$$
Hence (\ref{34}) holds.
Combining this case with the last one (i.e., $n=1$), it follows that
${\underline\kz}^{\rm DD}<4\,\kz^{\rm DD}$ for $\vz\in \big(\big(\sqrt{13}-1\big)/3,\, \sqrt{2}\,\big) $.
}\dexmp

\xmp{\bf\rf{cmf10}{Examples 7.7\,(5)}}\qd {\rm Let $E=\{1, 2, \cdots\}$, $a_i=1/i$ and $b_i=1$ for all $i\ge 1$.
Then $\lz^{\rm DD}=(3-\sqrt{5}\,)/2\approx 0.38$ and
 $\big(\lz^{\rm DD}\big)^{-1}\approx 2.618$.
We have
 $\mu_i=i!$, $\mu_i a_i=(i-1)!$ and $\mu_i b_i=i!$
 for all $i\ge 1$. Furthermore, we have
$$\big(\kz^{\rm DD}\big)^{-1}\!\!=\!\!\bigg(\bigg[\sum_{k=1}^1 \frac 1 {(k\! -\! 1)!}\bigg]^{-1}\!\! + \bigg[\sum_{k=4}^{\infty} \frac 1 {k!}
\bigg]^{-1}\bigg)\! \bigg[\sum_{\ell=1}^4 \ell !\bigg]^{-1}\!\!\!
=\!\frac 1 {33}\bigg[1 + \frac{3}{3 e\!-\!8}\bigg]\!\approx 0.6174.$$
And so $\kz^{\rm DD}\approx 1.62$.
With
$$\fz_i^-=\sum_{k=1}^i \frac{1}{\mu_k a_k}=\sum_{k=1}^i \frac{1}{(k-1)!}\qd
\text{and}\qd  \fz_k^+=\sum_{\ell=k}^\infty \frac{1}{\mu_{\ell} b_{\ell}}=\sum_{\ell=k}^\infty \frac{1}{\ell !},$$
we have
\begin{gather}
\dz_{m,\, \uz}^{-, \gz}=\sup_{i\le m}\bigg[\frac{1}{\gz} \uz!+\!\!\sum_{i\le \ell\le \uz-1} \ell!\bigg] \fz_i^-, \qqd
\dz_{n,\, \uz}^{+, \gz}=\sup_{i\ge n} \bigg[\frac{\gz-1}{\gz}\uz!+\!\!\sum_{\uz+1\le \ell \le i} \ell!\bigg] \fz_i^+.\nonumber\\
\frac{h_{\uz}^{-, \gz}}{f_m^{-, m}}
=\frac{1}{\sqrt{\fz_m^{-, \gz}}}\sum_{k=1}^{m-1}  (\fz_k^-)^{3/2} k!
+ \sum_{k=m}^{\uz-1}  \fz_k^- k! +\frac{1}{\gz} \fz_{\uz}^- \uz !,\nonumber\\
\frac{h_{\uz}^{+, \gz}}{f_n^{+, n}}
=\frac{1}{\sqrt{\fz_n^{+, \gz}}}\sum_{k=n+1}^\infty (\fz_k^+)^{3/2} k!
+ \sum_{k=\uz+1}^{n} \fz_k^+ k! +\frac{\gz-1}{\gz} \fz_{\uz}^+ \uz!.\nonumber
\end{gather}
For convenience, let $(L)$, $(R)$, $(M_-)$  and $(M_+)$ denote the last four quantities.

The candidates given by (\ref{32}) are $\uz=2,\,3$.
The case of $\uz=3$ is ruled out by (\ref{38}) and so we fix $\uz=2$. Then with $m=1,\,2$, $(m, n)$ satisfies (\ref{38}) for every $n$: $2\le n \le 17$. For the simplest choice
$m=n=\uz$, $\dz_{m,\, \uz}^{-, \gz}$ is attained at $i=1$
once $\gz\ge 2$, $\dz_{n,\, \uz}^{+, \gz}$ is attained at
$i=4$ whenever $\gz\ge 5/4$, and then the solution to the equation $\dz_{m,\, \uz}^{-, \gz}=\dz_{n,\, \uz}^{+, \gz}$ is $\gz\approx 3.2273$. Therefore,
$$\dz_{m,\, \uz}^{-, \gz}=\dz_{n,\, \uz}^{+, \gz}\approx 1.62.$$
A better choice is $(m, n)=(1, 5)$. Then
$\gz\approx 3.944$ and
$$4\times (R)=4\times (L)\approx 6.042,\qd (M_-)\approx 2.014, \qd (M_+)\approx 5.54.$$
This certainly implies (\ref{34}).}
\dexmp

\xmp{\bf\rf{cmf10}{Examples 7.7\,(8)}}\qd{\rm Let $E=\{1, 2, \cdots\}$,
$a_1=1$, $a_i=(i-1)^2$ for $i\ge 2$ and $b_i=i^2$ for $i\ge 1$.
Then $\lz^{\rm DD}=1/4=\big(4\,\kz^{\rm DD}\big)^{-1}$.
Once again, this example is dangerous.
Clearly,
$\mu_i=1$, $\mu_i a_i= a_i$ and $\mu_i b_i=b_i$ for all $i\ge 1$.
We have
$$\fz_i^-=1+\sum_{k=1}^{i-1} \frac{1}{k^2},\qqd \fz_i^+ =\sum_{\ell=i}^\infty \frac{1}{\ell^2}$$
and
\begin{gather}
\dz_{m,\, \uz}^{-, \gz}=\sup_{i\le m}\bigg[\frac{1}{\gz}+ \uz-i \bigg]\fz_i^-, \qqd
\dz_{n,\, \uz}^{+, \gz}=\sup_{i\ge n} \bigg[\frac{\gz-1}{\gz}+ i-\uz\bigg]\fz_i^+.\nonumber\\
\frac{h_{\uz}^{-, \gz}}{f_m^{-, m}}
=\frac{1}{\sqrt{\fz_m^-}}\sum_{k=1}^{m-1}  (\fz_k^-)^{3/2}
+ \sum_{k=m}^{\uz-1}  \fz_k^-  +\frac{1}{\gz} \fz_{\uz}^- ,\nonumber\\
\frac{h_{\uz}^{+, \gz}}{f_n^{+, n}}
=\frac{1}{\sqrt{\fz_n^+}}\sum_{k=n+1}^\infty (\fz_k^+)^{3/2}
+ \sum_{k=\uz+1}^{n} \fz_k^+ +\frac{\gz-1}{\gz} \fz_{\uz}^+.\nonumber
\end{gather}
As in the last example, we
use $(L)$, $(R)$, $(M_-)$  and $(M_+)$ to denote the last four quantities.
The only candidate by (\ref{32}) is $\uz=2$ which is fixed now. Then (\ref{38}) holds for all $m=1, 2$ and $n\ge 2$. The key for this example is that
$4\times (R)=4$, independent of $\uz$ and $\gz$. With $m=2$, the maximum of $(L)$
is achieved at $i=1$, it tends to 1 as $\gz\to\infty$. Since $m=\uz$, the term
$(M_-)$ is ignored. Besides, we have $(M_+)<4$ for all $n$: $2\le n\le 58$. Therefore,
(\ref{34}) holds.}
\dexmp

To conclude the paper, we make a remark on the generalization of the
results given here.

\rmk {\rm By using a known technique (cf. \rf{cmf05}{Section 6.7}), the variational
formula and its corollaries for the lower estimate of ${\lz}^{\text{\rm DD}}$ can be
extended to a more general setup (Poincar\'e-type inequalities). The upper estimate
is easier and was given in \rf{cmf10}{the remark above Corollary 8.3}.
}
\dermk

\nnd {\small School of Mathematical Sciences, Beijing Normal University,\newline
Laboratory of Mathematics and Complex Systems (Beijing Normal University),\newline
\text{\qqd} Ministry of Education,\newline
 Beijing 100875, The People's Republic of China.\newline
 E-mail: mfchen@bnu.edu.cn\newline
 Home page:
    http://math.bnu.edu.cn/\~{}chenmf/main$\_$eng.htm
}
\end{document}